\documentclass[12pt]{amsart}
\usepackage{amssymb,epsfig,amssymb,amsthm,epic,eepic,color,amsfonts}
\usepackage{hyperref}
\usepackage[all]{xy}
\usepackage{fontenc}  
\usepackage[utf8]{inputenc}

\newtheorem{thm}{Theorem}[section]

\newtheorem{prop}[thm]{Proposition}
\newtheorem{defn}[thm]{Definition}

\newtheorem{rien}[thm]{}
\numberwithin{equation}{section}

\newcommand{\be}{\begin{enumerate}}
\newcommand{\ee}{\end{enumerate}}
\newcommand{\bi}{\begin{itemize}}
\newcommand{\ei}{\end{itemize}}

\def\R{\mathbb{R}}

\def\ga{\gamma}

\def\be{\beta}

\def\De{\Delta}
\def\vp{\varphi}

\def\la{\lambda}

\def\ds{\displaystyle}
\def\p{\partial}

\def\nd{\noindent}
\def\bull{\hfill$\Box$\\}

\begin{document}
\vskip -1cm
\today
\bigskip

\begin{center}
{\sc A simple eversion of the 2-sphere with a unique quadruple point
}
\vspace{1cm}

Denis Sauvaget
\end{center}

\title{}
\author{}
\address{}
\email{denis.sauvaget@math.univ-paris13.fr}

\keywords{2-sphere eversion, immersion, regular homotopy, quadruple point}

\subjclass[2000]{57R19}

\begin{abstract} We give an example of an eversion of the 2-sphere in the Euclidean 3-space,
inspired by Morse theory, with a unique quadruple point. No homotopical tool is used.
\end{abstract}
\maketitle
%\tableofcontents
%\today
\thispagestyle{empty}
\vskip -2cm
\section{Introduction} 
We recall that an \emph{eversion} of the 2-sphere is a \emph{regular homotopy}, that is path of immersions 
of the unit 2-sphere into $\R^3$, starting from the Identity map of the 2-sphere and ending to a map $S^2\to S^2$
reversing the orientation. If the normal orientation points outwards at the beginning of the path it points inwards in the 
end.

\begin{center}
 \begin{figure}[h]
 \includegraphics[scale =.6]{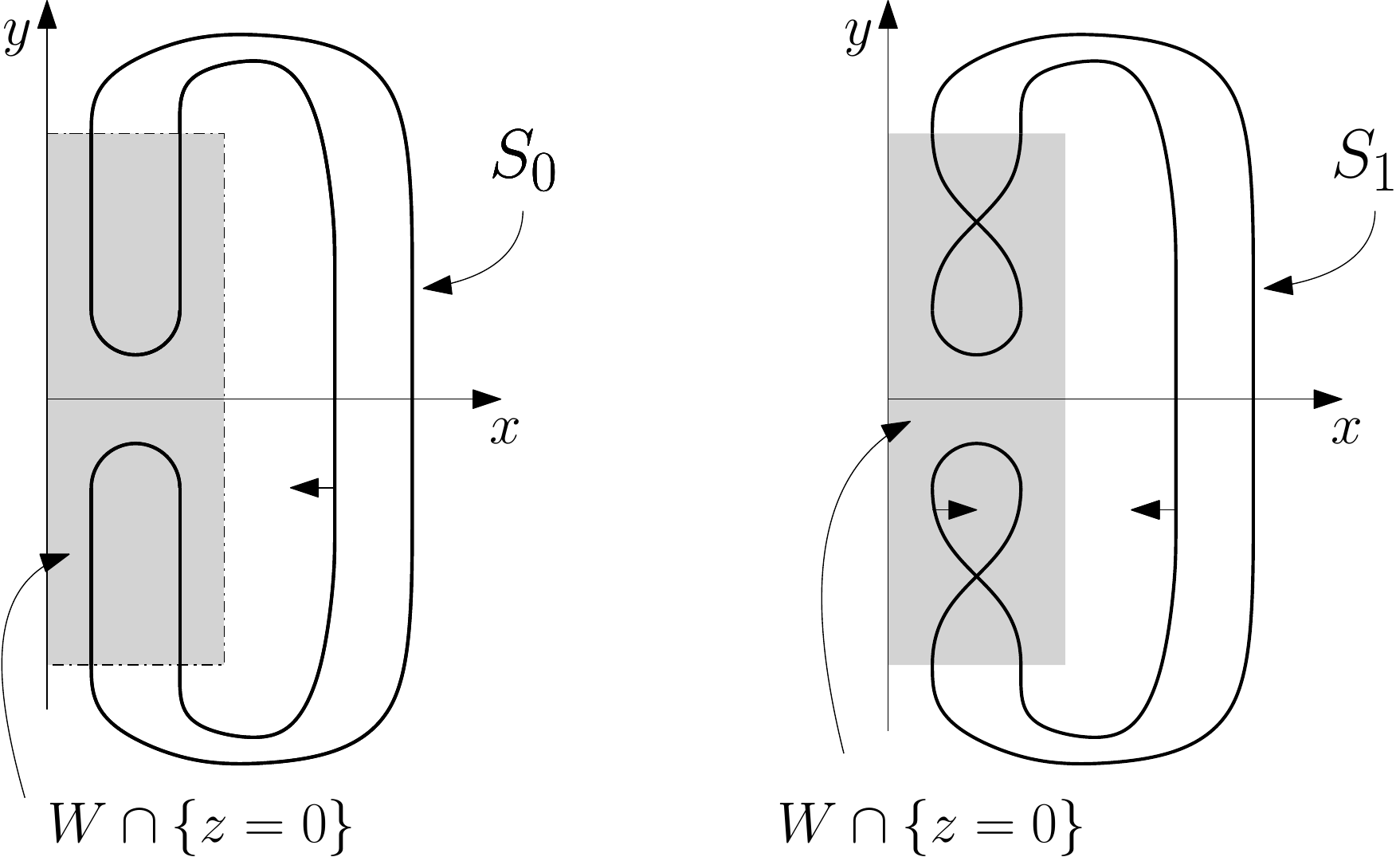}
 \caption{The trace of the homotopy from $S_0$ to $S_1$ on the plane $\{z=0\}$
 will be contained 
 in the  gray  rectangle.}
 \end{figure}\label{revol} 
 \end{center}

Since the unexpected outstanding result etablished by S. Smale \cite{smale}, 
many people gave examples of an eversion, including
H. Hopf \& N. Kuiper, A. Shapiro, M. Froissart, G. Francis \& B. Morin, F. Ap\'ery. 
The idea of Froissart-Morin, informally proved, circulated with the idea that there is a unique quadruple point
along their eversion. 
That allowed J. Hughes \cite{hughes} to prove that quadruple points are necessary in every path of immersions
realizing an eversion of the 2-sphere. In 1994, Ap\'ery \cite{apery} wrote a proof that there is
 a unique quadruple point in the 
Froissart-Morin eversion.\footnote{ It should be noted that the written proof is simplicial in nature 
and the smoothing is left to the reader.} More recently, A. Ch\'eritat \cite[2014]{cheritat} proposed another sphere
eversion that relies on the following property: namely, the space of  $C^1$ immersions of $[0,1]$ to the square 
with fixed extremities and  null winding number is contractible.

 In the present note, we construct an eversion where each step is controlled by a few planar figures.
First, we introduce two co-oriented surfaces  $S_0$ and $S_1$ of revolution about the $x$-axis in $\R^3$. 
They are represented by their respective co-oriented planar sections which are drawn in Figure \ref{revol}; $S_0$ is embedded and $S_1$  is immersed with a circle of double points. 

The surface $S_0$ is isotopic to the unit 2-sphere with its outward co-orientation. It is easy to check $S_1$ is regularly 
homotopic  to the unit 2-sphere with its inward co-orientation; indeed one cancel the circle of double points through 
the unique 3-ball whose angular boundary is contained in $S_1$.

We are going to prove the following proposition:

\begin{prop} \label{main}
The co-oriented surfaces $S_0$ and $S_1$ are regularly homotopic through a homotopy with a unique 
quadruple point.
\end{prop}
As a consequence, we have an eversion of the unit 2-sphere in $\R^3$ with a unique quadruple point.\\

\nd{\bf Acknowledgements.}
In 2022, the interest to present a new example of a sphere eversion (if it is really new) is not clear. 
Being not a specialist, I did it for convincing 
myself this move of the 2-sphere really exists.
 Thanks to my former thesis advisor Fran\c cois Laudenbach, I got information from 
%Professor 
Tony Phillips who promptly sent references to us, in particular about the necessity of quadruple points.
Then, 
%Professor 
Fran\c cois Ap\'ery had long conversations with Laudenbach about the history of that topic.
I am deeply grateful to all three of them.

\section{Support, Movies and  Routes}

\begin{rien}\label{support}
{\rm
The support of the regular homotopy from $S_0$ to $S_1$ that we are going to construct 
will be contained in a rectangular parallelepiped $W: =[0,2]\times [-6,6]\times[-6,6]$.
%where 
The units on each axis are such that, for $i=0,1$, the boundary of $W\cap S_i$ is made of two squares.
%and %each arcwise component of 
Moreover, $W\cap S_i$ is planar and parallel to $Oyz$ near its boundary .
}
\end{rien}
\begin{center}
 \begin{figure}[h]
 \includegraphics[scale =.6]{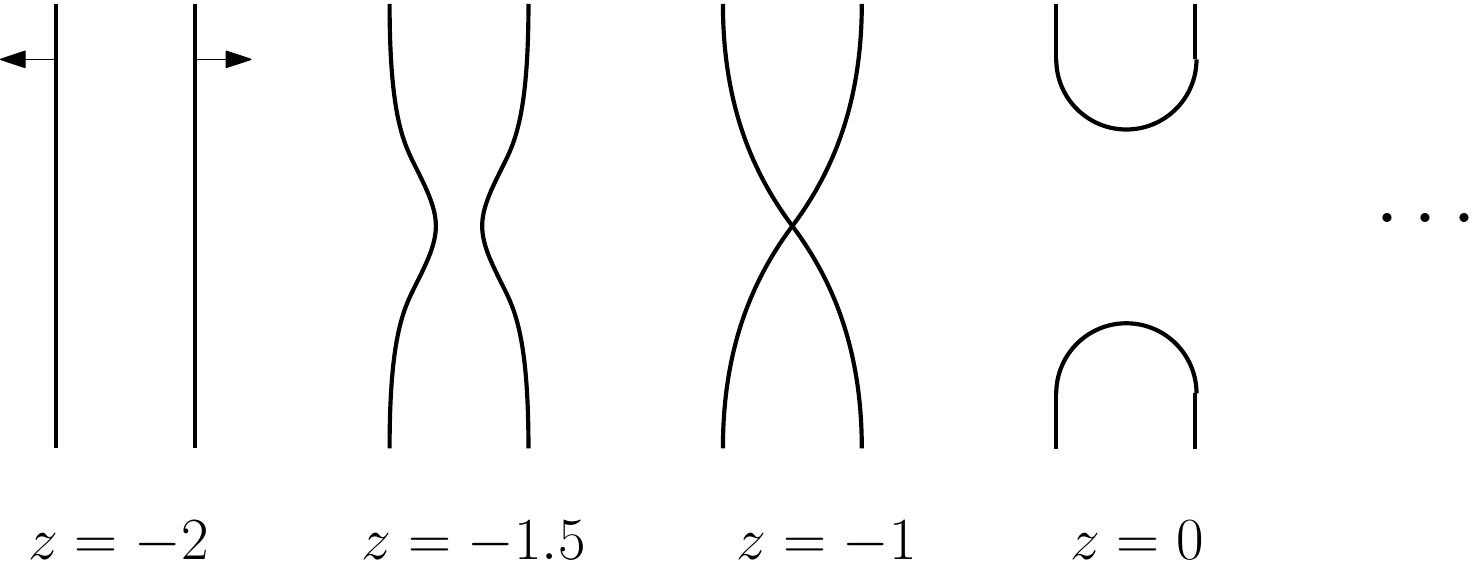}
 \caption{The movie is continued symmetrically from $z=0$ to $z=2$.}
 \end{figure}\label{C}
 \end{center}

\begin{defn} Given an immersed  surface $S\subset \R^3$ 
such that $W\cap S$ satisfies the same boundary condition as
 $W\cap S_0$, the \emph{movie} of $S$ is the family, parametrized by the Euclidean coordinate $z$ which plays the 
 role of time, of
the level sets of the \emph{height} function $z\vert_ {W\cap S}$.
\end{defn}
In Figure \ref{C}, we show a few moments of the movie of $S_0$. For a convenient choice of the unit of the 
$z$-axis, the critical values of the height function of $S_0$ are $\pm 1$ and for $z\leq -2$ or $z\geq 2$
the level set is made of two parallel segments. The lower saddle point has the sign $+$ meaning that $\p _z$
is an outward normal to $S_0$ at this point. The upper saddle point has the sign $-$.

\begin{center}
 \begin{figure}[h]
 \includegraphics[scale =.6]{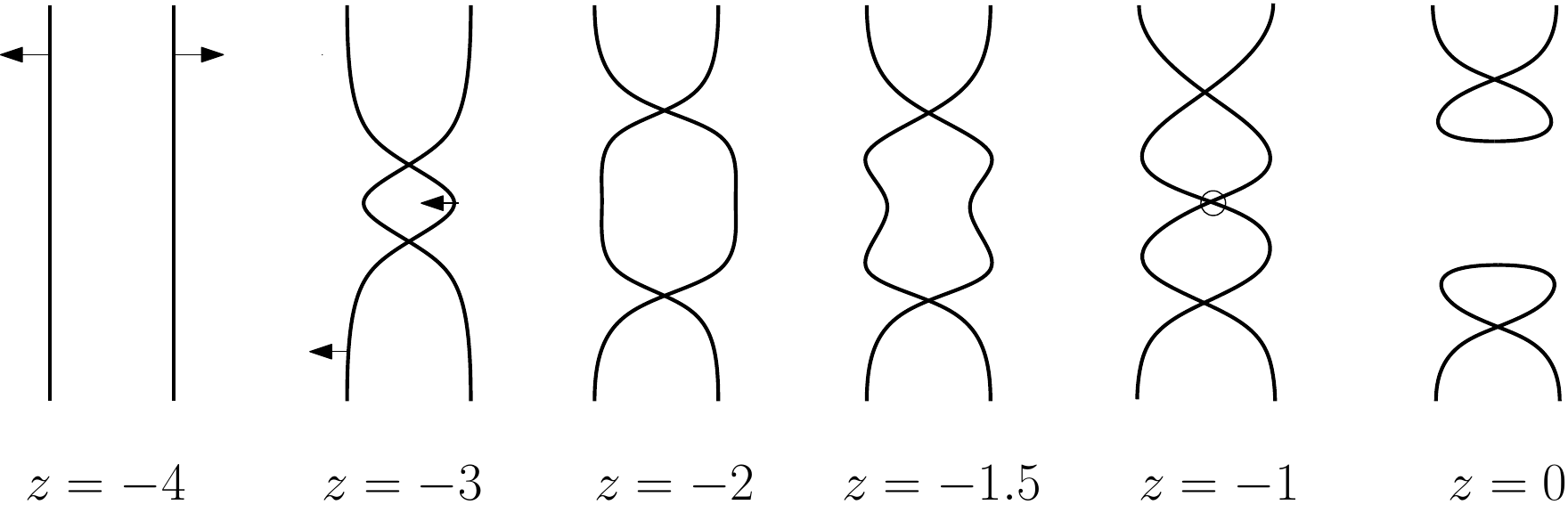}
 \caption{This movie is extended symmetrically from $z=0$ to $z=4$.
 Among the three double points at level $-1$, the saddle point is marked  with a circle.}\label{D}
 \end{figure}
 \end{center}

In Figure \ref{D}, we present the movie of $S_1$ from time $z= -4$ to $z=0$; it is extended by symmetry 
with respect to the horizontal plane $z=0$ to the interval $[0,4]$. The critical values of the height function
are still $z=\pm 1$. The lower saddle point has the sign $-$ and the upper saddle point has the sign $+$.

\begin{rien}
{\rm
The construction of the homotopy  $(S_t)_{t\in [0,1]}$  
deals only with $W$ (the \emph{support} of the desired regular homotopy) and is different in nature
for $z\leq 2$ and $z\geq 2$. % the latter domain will be filled up in Section \ref{section4}.
More precisely, the $z$ function, restricted to $S_t$,  
will have permanently two saddle points in $W$ and these critical points are moving in $W \cap\{z\leq 2\}$ only.
 For this part, it seems convenient to define the term of $route$. For short, in the remainder $S_t$ will denote what was
previously noted $S_t\cap W$.
}
\end{rien}

 \begin{defn} A route $+$ (resp. $-$)   
 is an embedded surface $R_{t}^+$ (resp. $R_{t}^-$),  contained in  $S_t$ with an octagonal  
 boundary and which satisfies the following conditions. 
  
 \begin{enumerate}
 \item the height function restricted to $S_t$ has a unique critical point in $R_t^\pm$ and that one is of type $\pm$. 
 \item The boundary $\p R_{t}^\pm$ is made of:
   \begin{itemize} 
  \item two opposite sides named the \emph{lower sides} %called the \emph{ends} 
 %attached with two distinct segments of a 
which lie at  some level $z_0$; % of the height; %$z\vert_{S_t}$, 
\item four vertical sides %segments
 between the levels %sets 
 $z_0$ and $z_0+2$,
 \item two sides of $R_t$ at level $z_0+2$ which are named the \emph{upper sides}. 
 \end{itemize} 
\item The projection to the plane  $Oxy$ of the interior of $R_t^\pm$ is a diffeomorphism and the  image of the octogon is a quadrilateral with cusps as vertices.
\end{enumerate}
 
 \end{defn} 
 
 This definition remembers the Morse model except the verticality conditions. 
 % and  that no non-degeneracy condition is imposed to the height function.

\begin{rien}\label{possible}
{\rm
Two possibilities  will be used   
for moving the routes and deforming $S_t$ in $W$:
\begin{itemize}
\item  
Vertically, $R_t^\pm$ is 
 moved by an isometry where  each point of $R_t^\pm$ remains on its own vertical line of $\R^3$.
\item  
Horizontally, $R_t^\pm$ is moved by an isotopy in which 
each level set of $R_t^\pm$ moves at a constant level and the vertical sides of $\p R_t^\pm$ %remain
 are kept vertical.
 \end{itemize}
 
  For extending such a horizontal isotopy to a \emph{regular homotopy} of $S_t\cap \{z\leq 2\}$ 
  %supported 
  in $W\cap\{z\leq 2\}$ one has
to remove some \emph{vertical rectangle}\footnote{ This means a surface, diffeomorphic to 
 a  rectangle, %(possibly non-planar) 
 with two vertical sides, two horizontal sides, and
 that is foliated by vertical segments.} 
swept out by one vertical 
side $A$  of $R_t^\pm$; at the same time, one glues some vertical rectangle swept out by 
another vertical side $A'\neq A$ of $R_t^\pm$
so that, at each time $t\in [0,1]$, $S_t$ is a proper immersed surface in $W$ with a fixed boundary. 

If $R_t^\pm$ moves vertically, $A$ is made of two horizontal arcs at the same level and $A'$ as well with $A'\cap A=\emptyset $\,.   Going down requiers some vertical room below $A$ %at the beginning 
for being swept out; and similarly
for moving up.

\begin{center}
 \begin{figure}[h]
 \includegraphics[scale =.7]{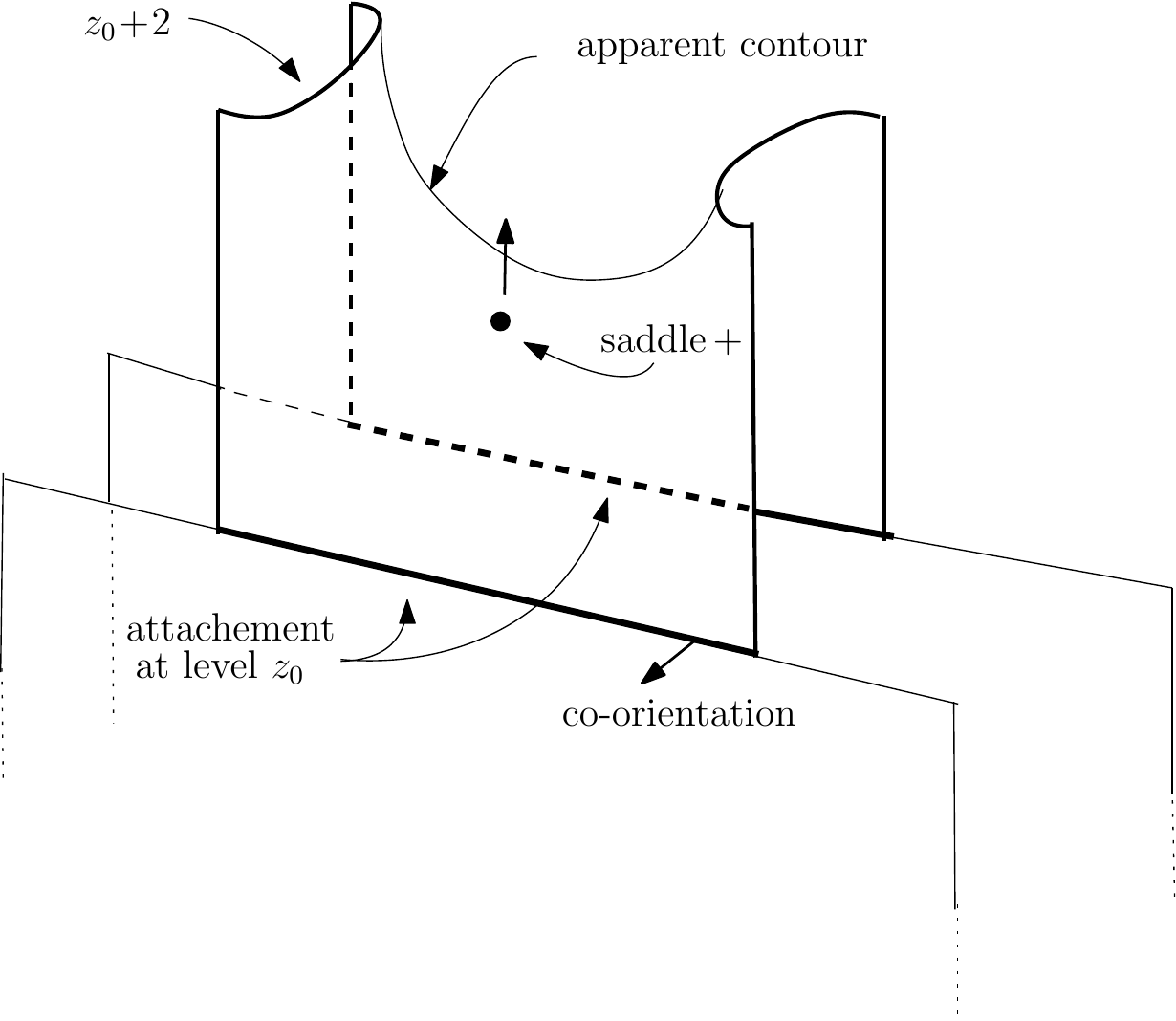}
 \caption{View of a route $+$.}\label{route}
 \end{figure}
 \end{center}
  }
\end{rien}
 
 \begin{center}
 \begin{figure}[h]
 \includegraphics[scale =.7]{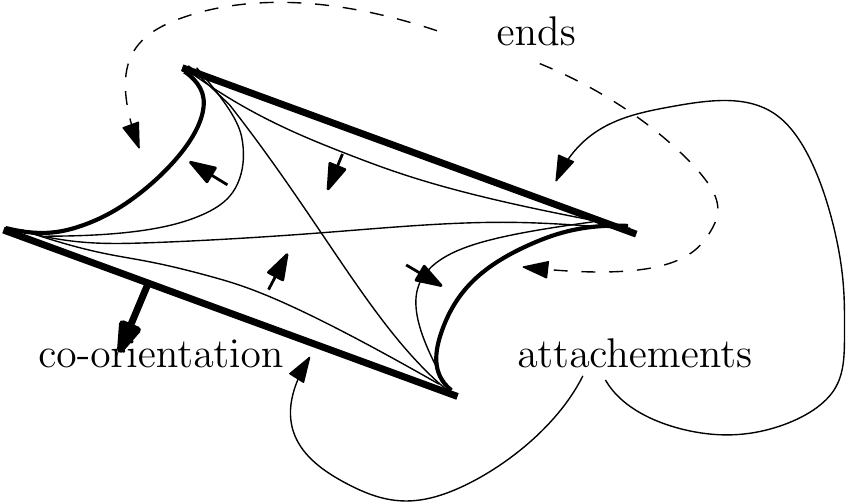}
 \caption{Projection to the
plane $Oxy$ of a route $+$. The thin lines
are the images of some level sets equipped with their ascending gradients.}
 \end{figure}\label{proj-route}\label{proj-route}
 \end{center}

The routes $-$ and $+$ have the same shape. Only the co-orientation is reversed; in other words, one turns Figure \ref{route}
 upside down.

\section{The  regular homotopy for $z\leq 2$} \label{section3}

A sequence of eleven  %ten
 figures is needed to describe successive steps of $S_t$ 
at times $t_1= 0< t_2< \cdots< t_{11}=1$  and mainly at the level $z=0$ (thick plain lines);
at every such a time, the projected routes to  the plane $z=0$ are drawn.
The route $+$ (resp. $-$) is drawn in red (resp. green). 
 When it seems to be useful, some other level sets of the routes will be drawn
 with dashed lines for understanding how the two routes, and hence the surface $S_{t_i}$, $i=1,\ldots,10$,
  are located in the domain  $z\in [-5,2]$.

\begin{rien}{\bf The first four figures.}
{\rm 
Figure \ref{t_1-to-t_4} (1) summaries Figure \ref{C}. Here, the two routes projects identically to the plane $z=0$;
the lower sides of the route $-$ (green) are exactly the upper sides of the route $+$ (red). That allows one to move 
the route $-$ horizontally  to the position drawn in Figure \ref{t_1-to-t_4} (2). The upper sides of $R^-_{t_2}$ are drawn 
with two dashed lines and lie at level $z=2$.

   \begin{center}
 \begin{figure}[h]
 \includegraphics[scale =.7]{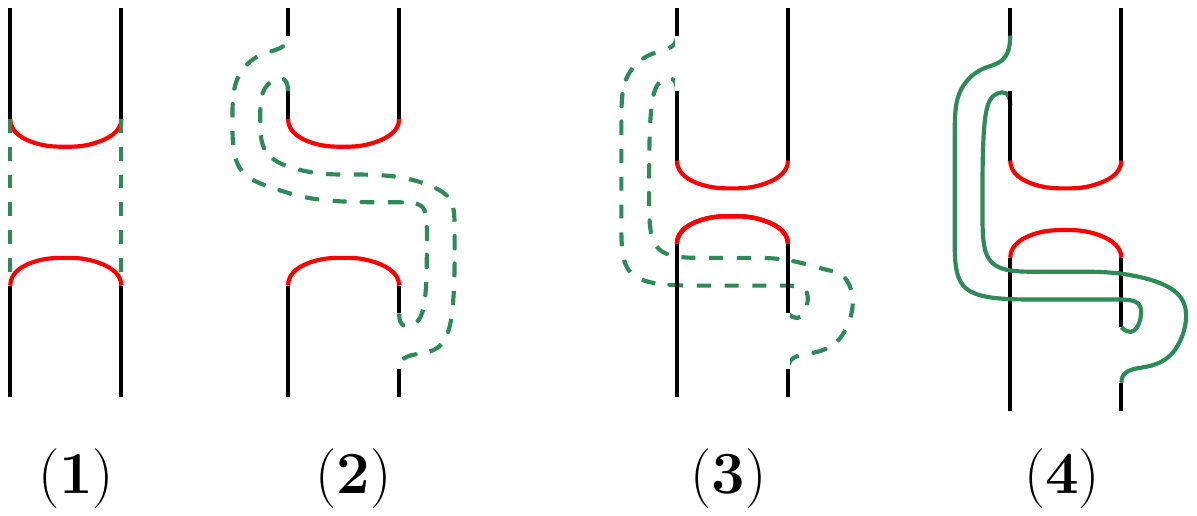}
 \caption{Here, the dashed lines are at level $z=2$.}\label{t_1-to-t_4}
 \end{figure}
 \end{center}
%}\end{rien}\end{document}
From Figure \ref{t_1-to-t_4} (2) to Figure \ref{t_1-to-t_4} (3), some horizontal isotopy is performed until
the projection of the green route avoids that of the red route.
Its movement during the next step, from (3) to (4), 
 is just a descending isotopy which is allowed since, at time $t_3$, the lower sides
of the green route are disjoint from the red route and have vertical rectangles below them.\footnote{ The projection  of
these rectangles to $Oxy$ are the small black dotted segments of Figure \ref{t_1-to-t_4} (3).}
The upper sides of $R^-_{t_4}$ lie at level 0. So, there are drawn with two plain lines.
}
\end{rien}

\begin{center}
 \begin{figure}[h]
 \includegraphics[scale =.7]{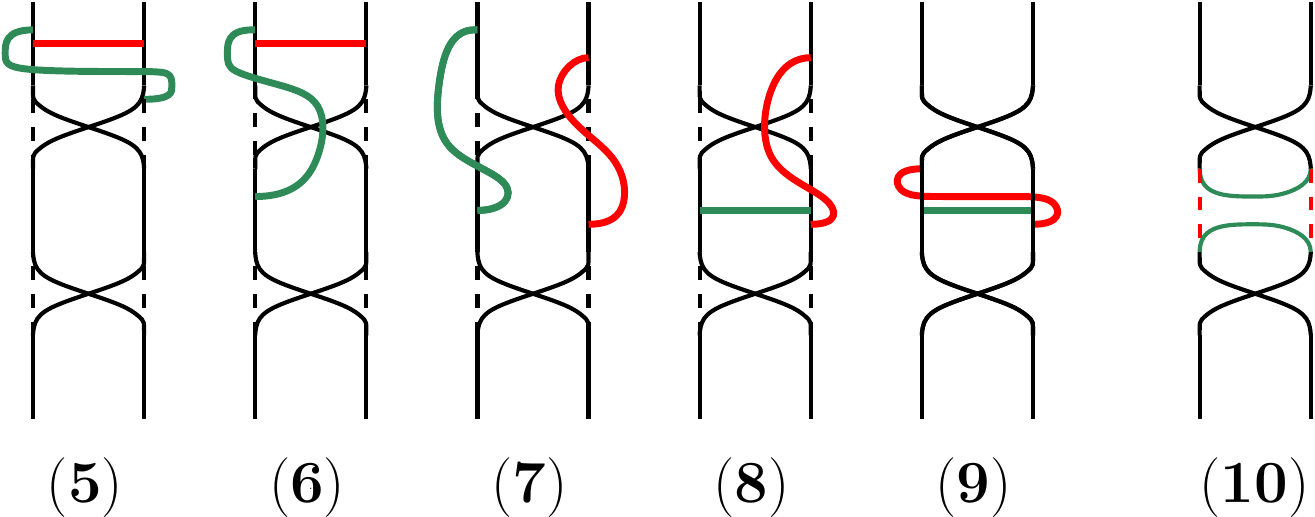}
 \caption{From $\bf (5)$ to $\bf (10)$ the dashed lines indicates the apparent contour of $S_t$ viewed from $z=\infty$.  
 In $\bf(11)$  
 the drawing is like Figure \ref{t_1-to-t_4}; the red dashed lines show the lower sides 
 of the red route, its upper sides coinciding with the green lines in $\{z=0\}$.}\label{t_5-to}
 \end{figure} 
 \end{center}

\begin{rien}{\bf The next seven figures.}
{\rm The two routes  remain disjoint and located  between the level sets $z=-2$ and $z=0$, except 
when $t\in [t_{10}, 1]$---see the end of the present subsection.
On Figure \ref{t_5-to}, for simplicity, 
the two routes are represented at each time $t_i$ %step
 by a thick colored plain line $L^\pm_{t_i}$, the sign recalling the one of the considered route. 
This line stands for the projection of the route 
$R^\pm_{t_i}$ to the plane 
$z=0$. The upper sides of $R^\pm_{t_i}$ are two lines essentially parallel to $L^\pm_{t_i}$ except near their 
extremities.

The regular homotopy $t\in [t_4,t_5]\mapsto S_t$ keeps the two routes pointwise fixed and consists of pushing
the finger from the left to the right on the left wall while the opposite movement is applied to 
the the right wall. 
At each time, the movie of $S_t$, $t\in [t_4,t_5]$,  is similar (that is, smoothly conjugate)  to the 
beginning of Figure \ref{D} (the first two steps), apart from the routes.
}
\end{rien}  
 
 \vskip -.8cm
  The steps from $t_5$ to $t_{10}$ need no special comments. Figures \ref{t_5-to} (5) to \ref{t_5-to} (10) 
  are quite explicit for describing each
movement. The last step, from (10) to (11), consists 
just of the series of three steps in Figure \ref{t_1-to-t_4} performed backwards, up to intertwining the colors. 
%Implicitly, 
The desired regular homotopy is now known in $W\cap\{z\leq 2\}$. \bull

For completing this homotopy
to the missing part $W\cap \{2\leq z\leq 6\}$ (see section \ref{section4}), it is useful to picture the two pathwise connected components of 
$S_t\cap \{z=2\}$. This is done in the next subsection.

\begin{rien} \label{3.3}{\bf The two components of $S_{t_i}$, $i=5, ...,11$, in level $z=2$.}
    \begin{center}
 \begin{figure}[h]
 \includegraphics[scale =.7]{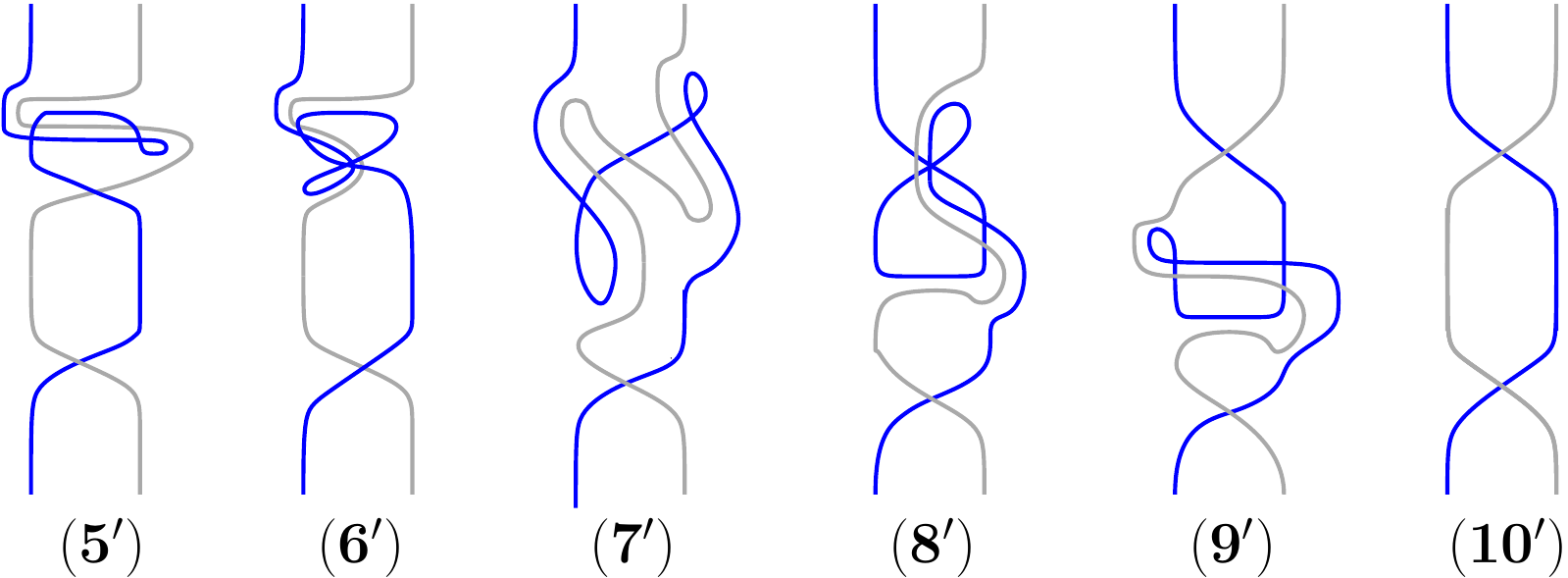}
 \caption{The two connected components, one in blue and the other in gray, 
 of $S_t\cap\{z=2\}$ at the corresponding times of Figure \ref{t_5-to}.}\label{2-comp}
 \end{figure}
 \end{center}
 
 {\rm %First, note that, at each time $t\geq t_4$, $S_t$ is vertical from $z=0$ up to $z=2$. So,
  The arcwise 
 connected components of
 $S_t\cap\{z= 2\}$ are easily deduced from the routes whose projection 
 are drawn on Figures \ref{t_1-to-t_4} and \ref{t_5-to}.
 %Implicitly, t
 The desired regular homotopy is now known in the domain $W\cap\{z\leq 2\}$. It will be completed 
 in the next section.
 }
 \end{rien}

\section{Completing $S_t$ in $\{2\leq z\leq 6\}$}\label{section4}

\begin{center}
\begin{figure}[h]
 \includegraphics[scale =.7]{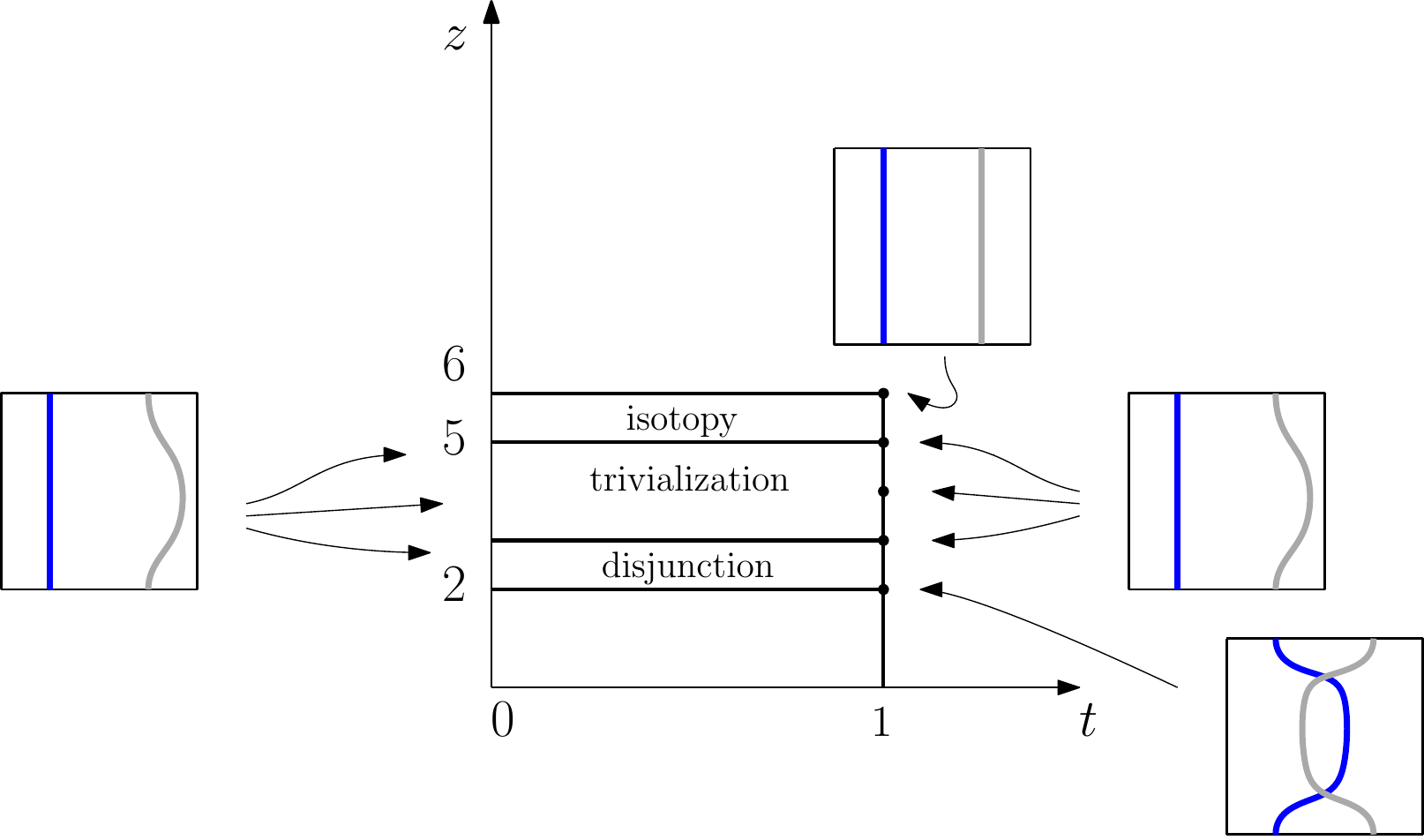}
 \caption{The boxes associated with some points in $(t,z)$-coordinates stand for $S_t\cap P_z\cap W$ equipped with 
 $(x,y)$-coordinates.}
 \label{plan-section4}
\end{figure}
 \end{center}

Let $P_z$ denote the horizontal plane $\R^2_{x,y}\times \{z\}$.
We explain in Figure \ref{plan-section4} how the movement is a function of time $t\in[0,1]$ 
depending on the interval  in which $z$ is located, either %$(-6, 2]$ or 
$[2,3]$ (\emph{disjunction zone}), $[3,5]$  (\emph{trivialization zone})
or $[5,6]$ (\emph{isotopy zone}). 
These different zones will be described in the remainder of the present section.
The quadruple point will appear in the disjunction zone.

The intersection $S_t\cap P_6\cap W$   is made of two fixed parallel segments for every $t\in [0,1]$; 
that is one of the boundary conditions
for gluing with the non-moving part of $S_t$. But the family $(S_t)_t$, that will be built in each zone,
does not exactly fulfill the right initial (resp. final) condition $S_0$ (resp. $S_1$) given in Figure \ref{revol}.
whose $z$-movies are given in Figures \ref{C} and \ref{D}; both are slightly bashed up (see the $(x,y)$-boxes 
in Figure \ref{plan-section4}.)
 Nevertheless, a banal
ambient isotopy, independent of $t$, rectifies this default. So, we neglect this phenomenon in the following discussion.

 \begin{rien} \label{disj}{\sc Disjunction zone.} {\rm This zone corresponds to $z\in [2,3]$.
 Let $b_t$ and $g_t$  ($b$ like blue and $g$ like
  grey---see Figure \ref{2-comp}) denote the
  two components at time $t$ of the curve $S_t\cap P_2\cap W$. \\
  
  \nd {\sc Claim.} {\it There exists
  an ambient  isotopy $(\vp^t_\la)_{\la\in[0,1]}$ of $P_2$,  supported in $P_2\cap int(W)$ and depending smoothly on $t$,
  which is viewed as an exernal parameter,
  such that $\vp^t_1(g_t)$ is disjoint from $b_{t'}$ for every $t'\in [0,1]$. Moreover,
  one may choose $(\vp^t_\la)$ such that the following holds:
  \begin{enumerate}
  \item The vector $\ds{\frac {d\vp^t_\la}{d\la}}$ points to the right of $\vp^t_\la(g_t)$ for every $(\la,t)\in [0,1]\times[0,1]$.
  \item For $\la=1$, the curve $\vp^t_1(g_t)$ is independent of $t$ and denoted by $\ga$.
  \item The domain of $P_2\cap W$ to the left of $\ga$ becomes convex after removing two strips
  $A^-$ and $A^+$, parallel to the $x$-axis, such that $g_t\cap A^\pm$ is independent of $t\in [0,1]$.\\
  \end{enumerate}
  }
  The third item imposed to $\ga$ will be used for the construction in the \emph{trivialization zone}.\\
  
   \nd{\sc Proof.} %One imposes this isotopy to lie in an autonomous flow whose infinitesimal generator $V^t$, depending on $t$ but not on $\la$,  points to the right of $g^t$ for every $t\in[0,1]$. 
  For $t=0$,  
  consider the one-dimensional foliation $\mathcal F_0$ parallel to the $x$-axis of the 2-dimensional domain 
  $R_0$ in $P_2\cap W$
  located to the right of the line  
  $g_0$; orient $\mathcal F_0$ as the $x$-axis.
  Choose an arbirary ambient isotopy $(\rho_t)_{t\in [0,1]}$ in
  $P_2$ supported in some compact set $C$ in $int(W)$ 
  such that $\rho_t(g_0)=g_t$. Then, we have an oriented foliation 
  $\mathcal F_t:= (\rho_t)_*(\mathcal F_0)$ of $R_t:= \rho_t(R_0)$. 
  We impose  the vector field $V^t$ to be tangent to $\mathcal F_t$ and point inwards $R_t$
  along $g_t$ for every $t$. It remains to choose the velocity
  of the desired autonomous flow.
  
  Consider the \emph{partial} boundary $\De$ of $R_0$ made of the complete boundary $\p R_0$ except the interior
  of the line $g_0$. There exists a collar $N$ of $\De$ in $R_0$ which fulfills the  following two conditions 
  for every $t\in[0,1]$:
  \begin{enumerate}
  \item[(i)]  $N\cap C=\emptyset$, that is, $N$ is kept pointwise fixed by $\rho_t$, and hence $N   \subset R_t$.
  \item[(ii)] $N$ is disjoint from $b_t$ for every $t\in [0,1]$. 
  \end{enumerate}
  Then $N$ is foliated by oriented segments parallel to the $x$-axis; 
  %starting from $\p N\smallsetminus \De$ and ending to $\De$; 
  this foliation is denoted by $\mathcal F_N$.
  One chooses a smooth path  $\ga$ in $N$, transverse to $\mathcal F_N$, joining the two end points of $g_0$
  and having the same germ as $g_0$ near the end points. 
  This arc $\ga$ can be chosen so that item (3) of the claim holds.

  %\begin{equation}\label{convex}\left\{begin{array}{l} \text{The domain of }P_2\cap W\text{ to the left of }\ga \text{ fulfills the following condition:}\\\text{the considered domain becomes convex after removing two strips } A^-\text{ and }A^+ \text{parallel to the }x\text{-axis, such  that }g_t\cap A^\pm \text{ is independent of } t\in [0,1]. \end{array}     \right.\end{equation}Two such strips can be easily found.
  
  Now, the vector field $V^0$ is chosen such that its flow, which reads 
  \begin{equation}
  \frac{d }{d\la}\vp^0_\la(x,y)=V^0\bigl(\vp^0_\la(x,y))\bigr),
   \end{equation}
   maps $g_0$ to $\ga$ in time 1.
   Then, if $V^t$ is the direct image of $V^0$ by $\rho_t: R_0\to R_t$, that is $V^t=\left(\rho_t\right)_*V^0$,
   its $\la$-flow $\vp^t_\la$ also maps $g_t$ to $\ga$
   in time 1. That completes the proof of the claim.
   ${}$\bull

The surface $S_t\cap \{2\leq z\leq 3\}$ has two components that are built as follows: one is vertical over $b_t$ 
 and the other  
 intersects the plane $P_z$ along the curve $\vp^t_{z-2}(g_t)$ 
 for every $t\in [0,1]$ and $z\in [2,3]$. Here, one uses the isotopy $\vp^t_\la$ from the claim.
 
 Since $g_t$ has no double point, a quadruple point only appears when $\vp_\la(g_t)$ meets a triple point of $b_t$. 
 By Subsection \ref{3.3},
 only   $b_{t_6}$ and $b_{t_9}$ have triple points. Moreover,
 the triple point of $b_{t_6}$ already lies on the left of $g_{t_6}$. Since the isotopy $\vp^{t_6}_\la$ constantly moves 
 $g_{t_6}$ to its right, the triple point of $b_{t_6}$ is never overlapped by $\vp^{t_6}_\la(g_{t_6})$. 
 By the same argument, the triple point
 of $b_{t_9}$ is overlapped exactly once by $\vp^{t_9}_\la(g_{t_9})$---see
   (9') in Figure \ref{2-comp}. So, there is exactly one quadruple point 
 of the family $(S_t)_t$ in the disjunction zone. Being careful in the other zones,  no quadruple points will be created 
 therein.\bull
  }
\end{rien}

\begin{rien} {\sc Trivialization zone.} {\rm This zone is defined by $z\in [3,5]$; it is divided in two parts.
Firstly, for $z\in [3,4]$ and $t\in [0,1]$,
one applies to the blue line $b_t$ an ambient isotopy  
$(\psi^t_\mu)_{\mu\in[0,1]}$ in the planar rectangle  $P_3\cap W$
 supported in the complement
of the fixed grey line $\ga=\vp^t_1(g_t)$ (see Claim from Subsection \ref{disj}) 
and ending in the \emph{normalized} position that is shown
in Figure \ref{vvv}. The parameter of this isotopy is $\mu$ and $t$ is an external parameter. 

The normalized blue line at time $t$ is contained in $P_4$ and denoted by $b'_t$. 
Its main property is to be immersed  and to have
at most two critical points of the  $y$-function restricted to $b'_t$.\footnote{ These two critical points are cancellable 
through immersions as shown by the $t$-movie either in  $[0,\frac{1}{10}]$ or in %the future to 
$[\frac{9}{10},1]$.}.

\begin{center}
\begin{figure}[h]
 \includegraphics[scale =.7]{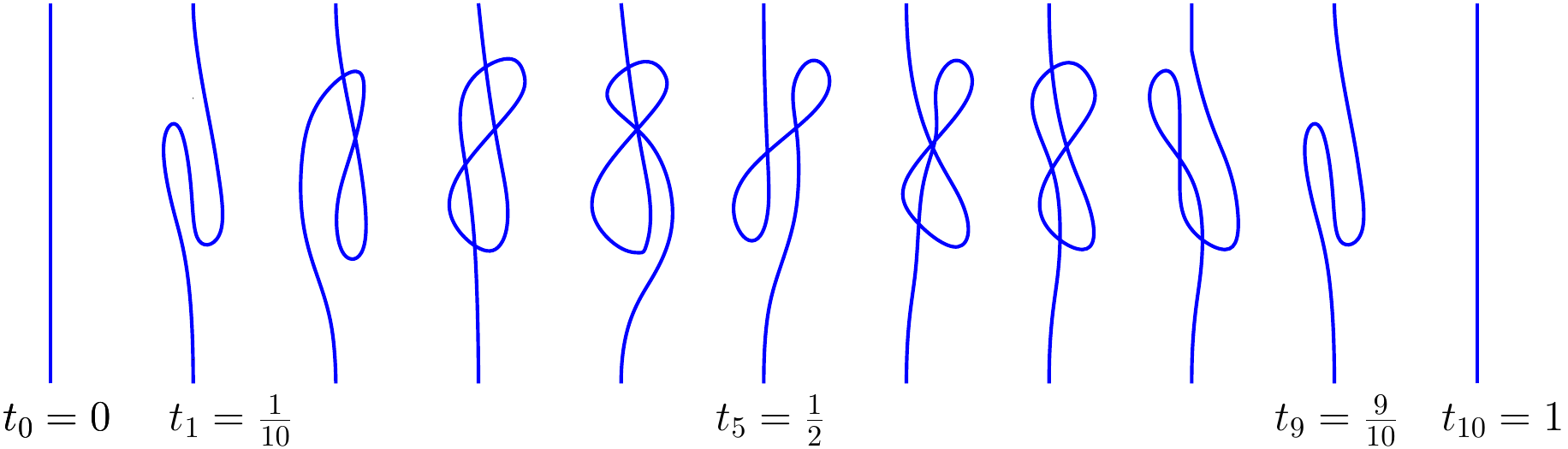}
 \caption{Here is a discrete $t$-movie in $P_4$ at times $t_i= \frac{i}{10}$, $i=0,1,\ldots,10$.}\label{vvv}
\end{figure}
 \end{center}

 In the domain $\{3\leq z\leq 4\}\cap W$, the surface $S_t$ has one fixed component which is vertical over $\ga$
 and one moving component whose level set at time $t$ is the line $\psi^t_{z-3}(b_t)\subset P_z$.
 The actual trivialization process starts from the normalized position in $\{z=4\}$ and consists of straightening $b'_t$
 smoothly in $t$. This will be performed in the next subsection.
 }
 \end{rien} 
 
 \begin{rien}{\sc Actual trivialization process.} \label{triv-process} {\rm One considers the family 
 $(b'_t)_{t\in [0,1]}$ of lines in the plane $P_4$ from Figure \ref{vvv}. Note that this plane is 
 canonically equipped with the $(x,y)$-coordinates  of the box $W$ which carries the family $S_t$
 of surfaces we are looking for.
 We choose a smooth parametrization of $b'_t$
 $$s\in [0,1]\mapsto \beta(s,t)=(x(s,t), y(s,t))\in P_4$$
  satisfying the
 following requirements:
 \begin{enumerate}
 \item For every $t\in[0,1]$, the map $s\mapsto \beta(s,t)$ is an immersion.
 \item The sign of the derivative $\p_sy$  changes exactly at $s= \frac13$ and $s=\frac 23$ for every $t\in
 \left[\frac{1}{10}, \frac{9}{10}\right]$. 
 \item $\beta(s,\frac{9}{10})=\beta(s,\frac{1}{10})$.
 \end{enumerate}
 In terms of these data, the trivialization is explicitely given by a barycentric combination formula that, 
 for every $z\in [4,5]$, reads
 \begin{equation}\label{bary}
 \left\{
 \begin{array}{lll}
 \hat\beta(z,s,t)&= (5-z) \beta(s,t)+(z-4)\beta(s,\frac{1}{10})\quad \text{if } t\in [\frac{1}{10},\frac{9}{10}]\\
                        &        \\
                        &= \beta(s,t)\quad \text{ if }t\in [0, \frac{1}{10}] \cup  [\frac{9}{10}, 1]\,.
  \end{array}    
  \right.                  
 \end{equation}  
 Of course, these two formulas coincide in their common  %definition 
 domains making their union continuous,
 but not smooth with respect to $t$. More precisely, the map $t\mapsto \hat \beta(-,-,t)$ is continuous with values in the space of $C^1$-functions of $(z,s)$, which is what we need.
 Nevertheless, the smoothing along $[4,5]\times \{\frac {i}{10}\}$, $i= 1 \text{ or }9$, is elementary, but useless.

 If $z\in[4,5]$ and $t\in [\frac{1}{10},\frac{9}{10}]$ are fixed, then $s\mapsto \hat\beta(z,s,t)\in P_z$ is an immersion into $P_z$.  Indeed, for $s\neq \frac 13, \frac 23$ the partial derivative $\p_s\hat\beta(z,s,t)$ does not vanish
 due to its $y$-component and, for $s=\frac 13  \text{ or } \frac 23 $, that is due to the $x$-component.
 In particular, for every fixed $t$ the map $(s,z)\in [0,1]\times [4,5]\mapsto (\hat \beta(f(z),s,t),z)$ parametrizes
 a proper  immersed
  connected component of $S_t$ in $W\cap\{z\in [4,5]\}$
  whatever the smooth function $f$. Here we use the convexity condition (3) in the claim of Subsection \ref{disj}.\footnote{
   In the removed strips $A^\pm$, %by the definition of \emph{quasi-convex}, 
   the barycentric combination
  (\ref{bary}) is \emph{trivial}, that is, independent of $z$.}
  There are also defaults of smoothness with respect to $z$ along the frontier between two consecutive zones.
  This will be answered in the next subsection.
 ${}$\bull
 }
 \end{rien}

 \begin{rien}{\sc Isotopy zone, end of the proof of Proposition \ref{main}.}
 {\rm 
 In the rectangle $P_5\cap W$, we see the following $t$-movie  made of the gluing, rescaling and smoothing of 
 the movie from 0 to $\frac{1}{10}$ and then the one from $\frac{9}{10}$ to 1. 
 This consists a one-parameter family 
 of proper curves in $P_5\cap W$ on which the $y$-coordinate has no critical points or a cancellable pair of critical points.
 Name $b''_t$ this blue curve at time $t$. For every $t\in [0,1]$, the path parametrized by $z\in [5,6]$
 is the cancelling path issued from $b''_t$. At the same time, one can perform the straigthening of the 
 corresponding grey curve.
 
 We are now going to answer the smoothness question raised in the end of Subsection \ref{triv-process}, which
 appears along $\{z=i\}, \text{ for every } i= 2, 3, 4, 5$. The way of reasoning is independent of $i$. Near $i$, the variable
  $z$ will be replaced with $f(z)$ where $f$ is a $C^\infty$ function that is increasing, coincides with $Id$ far from
  $i$, fulfills $f(i)=i$ and  all of whose derivatives vanish at $i$.

 As mentioned in the beginning of  Section \ref{section4}, this finishes the proof of Proposition \ref{main} 
 and completes the example of sphere eversion we would like to explain.
 
 }
 \end{rien}

 \end{document}